\documentclass[11pt]{article}
\usepackage{amssymb,amsmath,epsfig}

\title{Finite index subgroups of fully residually free groups}
\author{\textsf{Andrey Nikolaev} \and \textsf{Denis Serbin}}
\date{}
\date{\textsf{August 28, 2008}}

\newcommand{\type}[1]{{\mathrm{Type\!}\left({#1}\right)}}
\newcommand{\types}[1]{{\mathrm{Types\!}\left({#1}\right)}}
\newcommand{\paths}[1]{{\mathrm{Paths\!}\left({#1}\right)}}
\newcommand{\spaths}[1]{{\mathrm{SPaths\!}\left({#1}\right)}}
\newcommand{\comp}{\mathrm{Comp}}
\newcommand{\comm}[2]{{\mathrm{Comm}_{#1}(#2)}}
\newtheorem{theorem}{Theorem}

\newtheorem{example}{Example}
\newtheorem{lemma}{Lemma}
\newtheorem{cor}{Corollary}
\newtheorem{remark}{Remark}

\newtheorem{prop}{Proposition}

\begin{document}

\maketitle

\begin{abstract}
Using graph-theoretic techniques for f.g.
subgroups of $F^{\mathbb{Z}[t]}$ we provide a criterion for
a f.g. subgroup of a f.g. fully residually free group to be
of finite index. Moreover, we show that this criterion can
be checked effectively. Also we obtain an analogue of Greenberg-Stallings
Theorem for f.g. fully residually free groups, and prove that a f.g.
non-abelian subgroup of a f.g. fully residually free group
is of finite index in its commensurator.
\end{abstract}

\tableofcontents

\section{Introduction}
\label{sec:1}

Fully residually free (or {\it freely discriminated} \cite{BMR1}, or
{\it $\omega$-residually free} \cite{Rem}, or {\it limit} \cite{Sela1, Sela2})
groups have been extensively studied over the last ten years.
Although appeared first in 60's (see \cite{B}) this class of groups
drew much attention because of its connection with equations over
free groups. Recall that a group $G$ is called {\it fully residually
free} if for any finitely many non-trivial elements $g_1, \ldots, g_n
\in G$ there exists a homomorphism $\phi$ of $G$ into a free group $F$,
such that $\phi(g_i) \neq 1$ for $i \in [1,n]$. There are other
definitions of these groups more or less convenient depending on the
setting.

This class of groups can be studied from several viewpoints using many
different techniques. In particular, f.g. fully residually free groups are
fundamental groups of graphs of groups of a very particular type and their
structure can be described using Bass-Serre theory (see \cite{MR2,KM1,KM2}).
These groups are relatively hyperbolic with respect to their
maximal abelian subgroups (see \cite{Dah}), which provides another tool of
studying them from a geometric viewpoint. It is known that f.g. fully
residually free groups these groups act freely on $\mathbb{Z}^n$-trees
(see \cite{MR3}), etc.

Our study of fully residually free groups relies heavily on the fact
proved by Kharlampovich and Myasnikov (see \cite{KM2}) that every finitely
generated fully residually free group is embeddable into $F^{\mathbb{Z}[t]}$,
the free exponential group over the ring of integer polynomials $\mathbb{Z}[t]$.
This group was introduced by Lyndon (see \cite{L}) and he proved that
$F^{\mathbb{Z}[t]}$ (and hence its subgroups) is fully residually free. It
follows that one way to understand the properties of these groups is to study
finitely generated subgroups of $F^{\mathbb{Z}[t]}$.

A new technique to deal with $F^{\mathbb{Z}[t]}$ became available
recently when Myasnikov, Remeslennikov, and Serbin showed that
elements of this group can be viewed as reduced {\it infinite}
words in the generators of $F$ (see \cite{MRS}). It turned  out that
many  algorithmic problems for  finitely generated fully
residually free groups can be solved by the same methods as in the
standard free groups. Indeed, in \cite{MRS2} an analog of
Stallings' foldings (see \cite{S1, KM}) was introduced for an arbitrary finitely
generated subgroup of  $F^{\mathbb{Z}[t]}$, which allows one to
solve effectively the membership problem in $F^{\mathbb{Z}[t]}$,
as well as in an arbitrary finitely generated subgroup of it.
Next, in \cite{KMRS} this technique was further developed to obtain
the solution of many algorithmic problems. In particular, it was proved
that for a f.g. subgroups $G, H$ and $K$ of
$F^{\mathbb{Z}[t]}$ such that $H,K \leq G$ there are only finitely many
conjugacy classes of intersections $H^g \cap K$ in $G$. Moreover, one can
find a finite set of representatives of these classes effectively. This implies
that one can effectively decide whether two finitely generated
subgroups of $G$ are conjugate or not, and check if a given
finitely generated subgroup is malnormal in $G$. Needless to say that
all these results can be reformulated for f.g. fully residually free groups.

In the present paper we further develop these methods focusing on the problems
involving finite index subgroups. It is worth mentioning that non-abelian f.g.
subgroups of f.g. fully residually free groups have finite index in their
normalizers - this fact follows immediately from Theorem 7 \cite{KMRS} (see
also \cite{BH}) - but no criterion for detecting finite index subgroups was
known. In this paper we provide such a
criterion which can be checked effectively given a finite presentation of a
f.g. fully residually free group and a finite generating set of its subgroup.
This allows us to draw several corollaries including an analogue of
Greenberg-Stallings Theorem for free groups.

The authors are extremely grateful to Alexei G. Miasnikov for insightful discussions
and many helpful comments and suggestions.

\section{Preliminaries}
\label{sec:2}

Here we introduce basic definitions and notations which are to be
used throughout the whole paper. For more details see
\cite{MRS,MRS2}.

\subsection{Lyndon's free $\mathbb{Z}[t]$-group and infinite words}
\label{subs:2.1}

Let $F = F(X)$ be a free non-abelian group with basis $X$ and
$\mathbb{Z}[t]$ be a ring of polynomials with integer coefficients
in a variable $t$. In \cite{L} Lyndon introduced a {\em
$\mathbb{Z}[t]$-completion} $F^{\mathbb{Z}[t]}$ of $F$,  which is
called now the Lyndon's free $\mathbb{Z}[t]$-group.

It turns out that $F^{\mathbb{Z}[t]}$ can be described as a union
of a sequence  of  extensions of centralizers \cite{MR2}
\begin{equation}
\label{eq:2.1.1} F = G_0 < G_1 < \cdots < G_n < \cdots,
\end{equation}
where $G_{i+1}$ is obtained from $G_i$ by extension of all cyclic
centralizers in $G_i$ by a free abelian group of countable rank.

In \cite{MRS} it was shown that elements of $F^{\mathbb{Z}[t]}$
can be viewed as {\em infinite words} defined in the following
way. Let $A$ be a discretely ordered abelian group. By $1_A$ we
denote the minimal positive element of $A$.  Recall that if $a, b
\in A$ then the closed segment $[a,b]$ is defined as
$$[a,b] = \{x \in A \mid a \leq x \leq b \}.$$
Let $X = \{x_i \mid i \in I\}$ be a set. An {\em $A$-word} is a
function of the type
$$w: [1_A,\alpha_w] \to X^\pm,$$
where $\alpha_w \in A,\  \alpha_w \geqslant 0$. The element
$\alpha_w$ is called the {\em length} $|w|$ of $w$. By
$\varepsilon$ we denote the empty word. We say that $w$ is {\em
reduced} if $w(\alpha )\neq w(\alpha +1)^{-1}$ for any
$1\leq\alpha <\alpha _w.$ Then, as in a free group, one can
introduce a partial multiplication $\ast$, an inversion, a word
reduction etc., on the set of all $A$-words (infinite words)
$W(A,X)$. We write $u\circ v$ instead of $uv$ if $|uv|=|u|+|v|.$
All these definitions make it possible to develop infinite words
techniques, which provide a very convenient combinatorial tool
(for all the details we refer to \cite{MRS}).

It was proved in \cite{MRS} that $F^{\mathbb{Z}[t]}$ can be
canonically embedded into the set of reduced infinite words
$R(\mathbb{Z}[t],X)$, where $\mathbb{Z}[t]$, an additive group of
polynomials with integer coefficients, is viewed as an ordered
abelian group with respect to the standard lexicographic order
$\leqslant$ (that is, the order which compares the degrees of
polynomials first, and if the degrees are equal, compares the
coefficients  of corresponding terms starting with the terms of
highest degree). More precisely, the embedding of
$F^{\mathbb{Z}[t]}$ into $R(\mathbb{Z}[t],X)$ was  constructed by
induction, that is, all $G_i$ from the series (\ref{eq:2.1.1})
were embedded step by step in the following way. Suppose, the
embedding  of $G_i$ into $R(\mathbb{Z}[t],X)$ is already
constructed. Then, one chooses a {\em Lyndon's set} $U_i \subset
G_i$ (see \cite{MRS}) and the extension  of cenralizers of all
elements from $U_i$ produces $G_{i+1}$, which is now also
naturally embedded into $R(\mathbb{Z}[t],X)$.

The existence of an embedding of $F^{\mathbb{Z}[t]}$ into the set
of infinite words  implies automatically the fact that all
subgroups of $F^{\mathbb{Z}[t]}$ are also subsets of
$R(\mathbb{Z}[t],X)$, that is, their elements can be viewed as
infinite words. From now on we assume the embedding $\rho :
F^{\mathbb{Z}[t]} \longrightarrow R(\mathbb{Z}[t],X)$ to be fixed.
Moreover, for simplicity we identify $F^{\mathbb{Z}[t]}$ with its
image $\rho(F^{\mathbb{Z}[t]})$.

\subsection{Reduced forms for elements of $F^{\mathbb{Z}[t]}$}
\label{subs:2.2}

Following \cite{MRS} and \cite{MRS2} we introduce various normal
forms for elements in $F^{\mathbb{Z}[t]}$ in the following way.

We may assume that the set
$$U = \bigcup_i U_i$$
is well-ordered. Let
$$U_i = \{u_{i_1}, u_{i_2}, \ldots \} \subset G_i,$$
be enumeration of elements of $U_i$ in increasing order. Denote by
$I_i$ the set of indices $i_1, i_2, \ldots$ of elements from
$U_i$. Now $g \in G_{n+1} - G_n$ has the following representation
as a reduced infinite word:
\begin{equation}
\label{eq:2.2.1} g = g_1 \circ u_{n_1}^{\alpha_1} \circ g_2 \circ
\cdots \circ u_{n_l}^{\alpha_l} \circ g_{l+1},
\end{equation}
where $\ n_1, n_2, \ldots, n_l \in I_n,\ \ g_k \in G_n,\ k \in
[1,l+1],\ \ [g_k, u_{n_k}] \neq \varepsilon$ (or $g_{k}=\varepsilon$)
$[g_{k+1},u_{n_k}] \neq \varepsilon$ (or $g_{k+1}=\varepsilon$),
$k \in [1,l],$ $|\alpha_k|
>> 0,\ k \in [1,l]$ (recall that $\alpha >> 0$ if $\alpha \in \mathbb{Z}[t] -
\mathbb{Z}$). Representation (\ref{eq:2.2.1}) is called {\em
$U_n$-reduced} if the ordered $l$-tuple
$\{|\alpha_1|,|\alpha_2|,\ldots,|\alpha_l|\}$  is maximal with
respect to the left lexicographic order among all possible such
representations of $g$.

\begin{example}
{\rm
Suppose $u = xyx \in U \cap F(X)$ and $g \in G_1 - G_0$. If
$$g = u^{2t} \circ (y x) \circ u^{3t}$$
then there exists another representation of $g$
$$g = u^{2t-1} \circ (x y) \circ u^{3t+1}.$$
The corresponding $2$-tuples are $(2t,3t)$ and
$(2t-1,3t+1)$. In the former one maximization
of exponents of $u$ goes from left to right, while in
the latter one from right to left.}
\end{example}

From (\ref{eq:2.2.1}) one can obtain another representation of
$g$. Fix any $u$ from the list $u_{n_1}, u_{n_2}, \ldots,
u_{n_l}$. Then
\begin{equation}
\label{eq:2.2.4} g = h_1 \circ u^{\beta_1} \circ h_2 \circ \cdots
\circ u^{\beta_p} \circ h_{p+1},
\end{equation}
where $\beta_j = \alpha_{m_j}, m_j \in [1,l], j \in [1,p],\ h_1 =
g_1 \circ u_{n_1}^{\alpha_1} \circ \cdots \circ g_{m_1},\ h_{p+1}
= g_{m_p + 1} \circ \cdots \circ g_{l+1},\ h_k = g_{m_k + 1} \circ
\cdots \circ g_{m_{k+1}}, k \in [2,p]$. Representation
(\ref{eq:2.2.4}) is called a {\it $u$-representation} or a {\it
$u$-form} of $g$. In other words, to obtain a $u$-form one has to
"mark" in (\ref{eq:2.2.1}) only nonstandard exponents of $u$.
Representation (\ref{eq:2.2.4}) is called {\em $u$-reduced} if the
ordered $p$-tuple $\{|\beta_1|,|\beta_2|,\ldots,|\beta_p|\}$ is
maximal with respect to the left lexicographic order among all
possible $u$-forms of $g$.

Observe that if (\ref{eq:2.2.4}) is a $u$-form for $g$ and $g$ is
cyclically reduced then obviously
\begin{equation}
\label{eq:2.2.6} (h_1 \circ u^{\beta_1} \circ h_2 \circ \cdots
\circ u^{\beta_p} \circ h_{p+1}) \circ (h_1 \circ u^{\beta_1}
\circ h_2 \circ \cdots \circ u^{\beta_p} \circ h_{p+1})
\end{equation}
is a $u$-form for $g^2$. So, we call (\ref{eq:2.2.4}) {\it
cyclically $u$-reduced} if (\ref{eq:2.2.6}) is $u$-reduced.
\begin{lemma} \cite{MRS2}
\label{le:2.2.4} For any given $u$-reduced form of $g \in G_{n+1}
- G_n, u \in U_n$, there exists a cyclic permutation of $g$ such
that  its $u$-reduced form is cyclically $u$-reduced.
\end{lemma}
Let $g \in G_{n+1} - G_n$ have a $U_n$-reduced form
$$g = g_1 \circ u_{n_1}^{\alpha_1} \circ g_2 \circ \cdots \circ u_{n_l}^{\alpha_l} \circ g_{l+1},$$
where $\ u_{n_1}, u_{n_2}, \ldots, u_{n_l} \in U_n,\ \ g_k \in
G_n,\ k \in [1,l+1],\ \ [g_k,u_{n_k}] \neq \varepsilon$ (or $g_{k}=\varepsilon$),
$[g_{k+1},u_{n_k}] \neq \varepsilon$ (or $g_{k+1}=\varepsilon$), $|\alpha_k|
>> 0,\ k \in [1,l]$. Now, recursively one has a $U_{n-1}$-reduced form for
$g_i$
$$g_i = g(i)_1 \circ u_{m_1}^{\beta_{m_1}} \circ g(i)_2 \circ \cdots \circ u_{m_s}^{\beta_{m_s}}
\circ g(i)_s,$$ where $u_{m_1}, \ldots, u_{m_s} \in U_{n-1},\
|\beta_{m_k}| >> 0, k \in [1,s],\ g(i)_k \in G_{n-1}, k \in
[1,s+1]$ and one can get down to the free group $F$ with such a
decomposition of $g$, where step by step subwords between
nonstandard powers of elements from $U_i$ are presented as
$U_{i-1}$-forms, $i \in [1,n]$. Thus, from this decomposition one
can form the following series for $g$:
\begin{equation}
\label{eq:2.3.1} F < H_{0,1} < H_{0,2} < \cdots < H_{0,k(0)} <
H_{1,1} < \cdots < H_{1,k(1)} < \ldots
\end{equation}
$$\ldots < H_{n-1,k(n-1)} < H_{n,1} < \ldots < H_{n, k(n)},$$
where $H_{j,1}, \ldots, H_{j,k(j)}$ are subgroups of $G_{j+1}$,
which do not belong to $G_j$ and $H_{j,i}$ is obtained from
$H_{j,i-1}$ by a centralizer extension of a single element
$u_{j,i-1} \in H_{j,i-1} < G_j$. Element  $g$ belongs to $H_{n,
k(n)}$ and does not belong to the previous terms. Series
(\ref{eq:2.3.1}) is called an {\it extension series} for $g$.

\smallskip

Using the extension series above we can decompose $g$ in the
following way: $g \in H_{n, k(n)}$ has a $u_{n, k(n)}$-reduced
form
$$g = h_1 \circ u_{n, k(n)}^{\beta_1} \circ h_2 \circ \cdots \circ u_{n, k(n)}^{\beta_l} \circ h_{l+1},$$
where all $h_j, j \in [1,l+1]$ in their turn are $u_{n,
k(n)-1}$-reduced forms representing elements from $H_{n, k(n)-1}$.
This gives one a decomposition of $g$ related to its extension
series. We call this decomposition a {\it standard decomposition}
or a {\it standard representation} of $g$.

\smallskip

Observe that for any $g \in F^{\mathbb{Z}[t]}$, its standard
decomposition can be viewed as a finite product $b_1 b_2 \cdots
b_m$, where
$$b_i \in B = \{X \cup X^{-1}\} \cup \{ u^\alpha \mid u \in U, \alpha \in \mathbb{Z}[t] - \mathbb{Z}\}.$$
We denote this product by $\pi(g)$ so we have
$$\pi(g) = \pi(h_1)\ u_{n,k(n)}^{\beta_1}\ \pi(h_2)\ \cdots\ u_{n, k(n)}^{\beta_l}\  \pi(h_{l+1}),$$
where $\pi(h_i)$ is a finite product in the alphabet $B$
corresponding to $h_i$, and from now on, by a standard
decomposition of an element $g$ we understand not the
representation of $g$ as a reduced infinite word but the finite
product $\pi(g)$.

\smallskip

By $U(g)$ we denote a finite subset of $U$ such that if $\pi(g)$
contains a letter $b_i \in B$ such that $b_i = u^\alpha$ then $u
\in U(g)$. Observe that $U(g)$ is ordered with an order induced
from $U$, so we have
$$U(g) = \{u_1,\ldots,u_m\},$$
where $u_i < u_j$ if $i < j$ and $u_m = u_{n,k(n)}$. By
$\max\{U(g)\}$ we denote the maximal element of $U(g)$.

\smallskip

If $u \in U(g)$ then by $deg_u(g)$ we denote the maximal degree of
infinite exponents of $u$, which appear in $\pi(g)$.

\smallskip

It is easy to see that in general $\pi(g_1 \circ g_2) \neq
\pi(g_1) \pi(g_2)$ and $\pi(g \circ g) = \pi(g) \pi(g)$ if and
only if the $u$-reduced form of $g$ is cyclically $u$-reduced,
where $u = \max\{U(g)\}$.

From the definition of a Lyndon's set and the results of
\cite{MRS} it follows that if $R \subset G_n$ is a Lyndon's set
then a set $R'$ obtained from $R$ by cyclic decompositions of its
elements is also a Lyndon's set. Thus, by Lemma \ref{le:2.2.4} we
can assume a $w$-reduced form of any $u \in U_n$ to be cyclically
$w$-reduced, where $w = \max\{U(u)\}$. Hence, we can assume
$$\pi(u \circ u) = \pi(u) \pi(u)$$
for any $u \in U$.

\subsection{Embedding theorems}

There are three results which play an important role in this
paper. The first embedding theorem is due to Kharlampovich and
Myasnikov.

\begin{theorem} [The first embedding theorem (\cite{KM2})]
\label{th:embKM} Given a finite presentation  of a finitely
generated fully residually free group $G$ one can effectively
construct an embedding $\phi:G \rightarrow F^{\mathbb{Z}[t]}$ (by
specifying the images of the generators of $G$).
\end{theorem}

Combining Theorem \ref{th:embKM} with the result on the
representation of $F^{\mathbb{Z}[t]}$ as a union of a sequence of
extensions of centralizers one can get the following theorem.

\begin{theorem} [The second embedding theorem]
\label{th:embKMR} Given a finite presentation  of a finitely
generated fully residually free group $G$ one can effectively
construct a finite sequence of extension of centralizers
 $$F < G_1 < \ldots < G_n,$$
where $G_{i+1}$ is an extension of the  centralizer of some
element $u_i \in G_i$ by an infinite cyclic group $\mathbb{Z}$,
and an embedding $\psi^\ast:G \rightarrow G_n$ (by specifying the
images of the generators of $G$).
\end{theorem}

Combining Theorem \ref{th:embKM} with the result on the effective
embedding of $F^{\mathbb{Z}[t]}$ into $R(\mathbb{Z}[t],X)$
obtained in \cite{MRS} one can get the following theorem.

\begin{theorem} [The third embedding theorem]
\label{th:embMR} Given a finite presentation  of a finitely
generated fully residually free group $G$ one can effectively
construct an embedding $\psi : G \rightarrow R(\mathbb{Z}[t],X)$
(by specifying the images of the generators of $G$).
\end{theorem}

\subsection{Graphs labeled by infinite $\mathbb{Z}[t]$-words}
\label{subs:2.3}

By an {\em $(\mathbb{Z}[t],X)$-labeled directed graph} ({\em
$(\mathbb{Z}[t],X)$-graph}) $\Gamma$ we understand
   a combinatorial graph $\Gamma$ where
every edge has a direction  and is labeled either by a letter from
$X$ or by  an infinite word $u^\alpha \in F^{\mathbb{Z}[t]}, u \in
U, \alpha \in \mathbb{Z}[t], \alpha > 0$, denoted $\mu(e)$.

For each edge $e$ of $\Gamma$ we denote the origin of $e$ by
$o(e)$ and the terminus of $e$ by $t(e)$.

For each edge $e$ of $(\mathbb{Z}[t],X)$-graph we can introduce a
formal inverse $e^{-1}$ of $e$ with the label $\mu(e)^{-1}$ and
the endpoints defined as $o(e^{-1}) = t(e), t(e^{-1}) = o(e)$,
that is, the direction of $e^{-1}$ is reversed with respect to the
direction of $e$. For the new edges $e^{-1}$ we set $(e^{-1})^{-1}
= e$. The new graph, endowed with this additional structure we
denote by $\widehat{\Gamma}$. Usually we will abuse the notation
by disregarding the difference between $\Gamma$ and
$\widehat{\Gamma}$.

A {\em path} $p$ in $\Gamma$ is a sequence of edges $p = e_1
\cdots e_k$, where each $e_i$ is an edge of $\Gamma$ and the
origin of each $e_i$ is the terminus of $e_{i-1}$. Observe that
$\mu(p) = \mu(e_1) \dots \mu(e_k)$ is a word in the alphabet $\{X
\cup X^{-1}\} \cup \{ u^\alpha \mid u \in U, \alpha \in
\mathbb{Z}[t] \}$ and we denote by $\overline{\mu(p)}$ a reduced
infinite word $\mu(e_1) \ast \cdots \ast \mu(e_k)$ (this product
is always defined).

A path $p = e_1 \cdots e_k$ in $\Gamma$ is called {\em reduced} if
$e_i \neq e_{i+1}^{-1}$ for all $i \in [1,k-1]$.

A path $p = e_1 \cdots e_k$ in $\Gamma$ is called {\em label
reduced} if
\begin{enumerate}
\item[1)] $p$ is reduced; \item[2)] if $e_{k_1} \cdots e_{k_2},\
k_1 \leq k_2$ is a subpath of $p$ such that  $\mu(e_i) =
u^{\alpha_i}, u \in U,  \alpha_i \in \mathbb{Z}[t],\ i \in
[k_1,k_2]$ and $\mu(e_{k_1 - 1}) \neq u^\beta, \mu(e_{k_2 + 1})
\neq u^\beta$ for any $\beta \in \mathbb{Z}[t]$, provided $k_1 -
1, k_2 + 1 \in [1,k]$, then $\alpha = \alpha_{k_1} + \cdots +
\alpha_{k_2} \neq 0$ and $\mu(e_{k_1 - 1}) \ast u^\alpha =
\mu(e_{k_1-1}) \circ u^\alpha, u^\alpha \ast \mu(e_{k_2 + 1}) =
u^\alpha \circ \mu(e_{k_2 + 1})$.
\end{enumerate}

Let $\Gamma$ be a $(\mathbb{Z}[t],X)$-graph and $u \in U$ be
fixed. Vertices $v_1, v_2 \in V(\Gamma)$ are called {\em
$u$-equivalent} (denoted $v_1 \sim_u v_2$) if there exists a path
$p = e_1 \cdots e_k$ in $\Gamma$ such that $o(e_1) = v_1, t(e_k) =
v_2$ and $\mu(e_i) = u^\alpha_i, \alpha_i \in \mathbb{Z}[t], i \in
[1,k]$. $\sim_u$ is an equivalence relation on vertices of
$\Gamma$, so if $\Gamma$ is finite then all its vertices can be
divided into a finite number of pairwise disjoint equivalence
classes. Suppose, $v \in V(\Gamma)$ is fixed. One can take the
subgraph of $\Gamma$ spanned by all the vertices which are
$u$-equivalent to $v$ and remove from it all edges with labels not
equal to $u^\alpha, \alpha \in \mathbb{Z}[t]$. We denote the
resulting subgraph of $\Gamma$  by $Comp_u(v)$ and call a {\em
$u$-component of $v$}. If $v \in V(\Gamma), v_0 \in V(Comp_u(v))$
then one can define a set
$$H_u(v_0) = \{\overline{\mu(p)} \mid p\ {\rm is\ a\ reduced\ path\
in}\ Comp_u(v)\ {\rm from}\ v_0\ {\rm to}\ v_0 \}.$$
\begin{lemma} \cite{MRS2}
\label{le:3.3.1} Let $\Gamma$ be a $(\mathbb{Z}[t],X)$-graph and
$v \in V(\Gamma), v_0 \in V(Comp_u(v))$.  Then
\begin{enumerate}
\item $H_u(v_0)$ is a subgroup of $R(\mathbb{Z}[t],X)$; \item
$H_u(v_0)$ is isomorphic to a subgroup of $\mathbb{Z}[t]$; \item
if $Comp_u(v)$ is a finite graph, then $H_u(v_0)$ is finitely
generated; \item if $v_1 \in V(Comp_u(v))$ then $H_u(v_0) \simeq
H_u(v_1)$.
\end{enumerate}
\end{lemma}

Following \cite{MRS2} one can introduce operations on
$u$-components which are called {\em $u$-foldings}.  One of the
most important properties of $u$-foldings is that they do not
change subgroups associated with $u$-components.
\begin{lemma} \cite{MRS2}
\label{le:3.4.4} Let $\Gamma$ be a $(\mathbb{Z}[t],X)$-graph, $v
\in V(\Gamma)$ and  $C = Comp_u(v)$ be finite. Then there exist a
$(\mathbb{Z}[t],X)$-graph $\Delta$ obtained from $\Gamma$ by
finitely many $u$-foldings such that $v' \in V(\Delta)$
corresponds to $v$ and $C' = Comp_u(v')$ consists of a simple
positively oriented path $P_{C'}$, and some edges that are not in
$P_{C'}$ connecting some pairs of  vertices in $P_{C'}$.
\end{lemma}
$C'$ in Lemma \ref{le:3.4.4} is called a {\it reduced}
$u$-component. Since  $P_{C'}$ is a simple path there exists a
vertex $z_{C'} \in V(P_{C'})$ which is an origin of only one
positive edge in $P_{C'}$. $z_{C'}$ is called a {\it base-point}
of $C'$.

It turns out that any finite reduced $u$-component $C$ in a
$(\mathbb{Z}[t],X)$-graph is characterized completely by the pair
$(P_C,H_u(z_C))$ in the following sense. For any reduced path $p$
in $C$ there exists a unique reduced subpath $q$ (denoted $q =
[p]$) of $P_C$ with the same endpoints as $p$, such that
$\overline{\mu(p)} \ast \overline{\mu(q)}^{-1} \in H_u(z_C)$.
Moreover, let $P_C = f_1 \cdots f_m$, where $o(f_1) = z_C, v_0 =
z_C, v_i = t(f_i), i \in [1,m]$ and let $p_0, p_1, \ldots, p_m$ be
reduced subpaths of $P_C$ such that $o(p_i) = z_C, t(p_i) = v_i, i
\in [0,m]$. The set of paths $p_0, p_1, \ldots, p_m$ is called a
{\em set of path representatives associated with $C$} (denoted by
$Rep(C)$).
\begin{lemma} \cite{MRS2}
\label{le:3.4.5}
Let $C$ be a finite reduced $u$-component in a
$(\mathbb{Z}[t],X)$-graph $\Gamma$,  $v \in V(C)$ and let $\alpha
\in \mathbb{Z}[t]$. If $\overline{\mu(p_i)} \ast
\overline{\mu(p_j)}^{-1} \notin H_u(z_C)$ for any $p_i, p_j \in
Rep(C), i \neq j$ then either there exists a unique reduced path
$p$ in $P_C$ such that $o(p) = v$ and $u^\alpha \in
\overline{\mu(p)} \ast H_u(z_C)$ or there exists no path $q$ in
$C$ with this property.
\end{lemma}
If $C$ is reduced and $Rep(C)$ satisfies the condition from Lemma
\ref{le:3.4.5} then we call  $C$ {\em a $u$-folded $u$-component}.

\subsection{Languages associated with $(\mathbb{Z}[t],X)$-graphs}
\label{subs:2.4}

Let $\Gamma$ be a $(\mathbb{Z}[t],X)$-graph and let $v$ be a
vertex of $\Gamma$.  We define the language of $\Gamma$ with
respect to $v$ as
$$L(\Gamma,v) = \{\overline{\mu(p)} |\ p~ {\rm is~ a~ reduced~ path~ in}~
\Gamma~ {\rm from}~ v ~{\rm to}~ v\}.$$

\begin{lemma} \cite{MRS2}
\label{le:4.1.1} Let $\Gamma$ be a finite
$(\mathbb{Z}[t],X)$-graph and let $v \in V(\Gamma)$.  Then
$L(\Gamma,v)$ is a subgroup of $F^{\mathbb{Z}[t]}$.
\end{lemma}

Let $\Gamma$ be a $(\mathbb{Z}[t],X)$-graph and $p = e_1 \cdots
e_k$ be a reduced path in $\Gamma$. Let $g \in G_{n+1} - G_n$ and
let
$$\pi(g) = \pi(g_1) u^{\beta_1} \pi(g_2) \cdots u^{\beta_l}
\pi(g_{l+1}),$$ be the standard decomposition of $g$, where $u =
\max\{U(g)\}$. We write
$$\mu(p) = \pi(g)$$
if $p$ can be subdivided into subpaths
$$p = p_1 d_1 p_2 \cdots d_l p_{l+1},$$
where each $d_i$ is a path in some $u$-component of $\Gamma$ so that
$\overline{\mu(d_i)} = u^{\beta_i}$, and each $p_i$
is a path in $\Gamma$ which does not contain edges labeled by
$u^\alpha,\ \alpha \in \mathbb{Z}[t]$ so that the equality
$\mu(p_i) = \pi(g_i)$ is defined inductively in the same way. Observe
that if $g = x_1 \cdots x_r \in F(X)$ then $\mu(p) = \pi(g)$ if and only if
$k = r$ and $\mu(e_i) = x_i,\ i \in [1,r]$.

\smallskip

Let $\Gamma$ be a finite $(\mathbb{Z}[t],X)$-graph. Since $\Gamma$
is finite, the set of elements $u \in U$ such that there exists an
edge $e$ in $\Gamma$ labeled by $u^\alpha, \alpha \in
\mathbb{Z}[t]$ is finite and ordered with the  order induced from
$U$. Thus one can associate with $\Gamma$ an ordered set
$U(\Gamma) = \{u_1,\ldots,u_N\}, N > 0, u_i \in U, u_i < u_j$ for
$i < j$.

\smallskip

Let $u_i \in U(\Gamma)$ be fixed and $\Gamma(i)$ be a subgraph of
$\Gamma$ which consists only of edges  $e \in E(\Gamma)$ such that
either $\mu(e) = x \in X^\pm$ or $\mu(e) = u_j^\alpha, \alpha \in
\mathbb{Z}[t], j \leq i$. $\Gamma(i)$ is called an {\em $i$-level
graph of $\Gamma$} (by $0$-level graph we understand a subgraph of
$\Gamma$ which consists only of edges with labels from $X$) and
the {\em level} (denoted $l(\Gamma)$) of $\Gamma$ is the minimal
$n \in \mathbb{N}$ such that $\Gamma = \Gamma(n)$. Observe that
$\Gamma(i)$ may not be connected for some $i < l(\Gamma)$, but
still one  can apply to $\Gamma(i)$ partial and $u$-foldings, $u
\in U(\Gamma)$.

\smallskip

A finite connected $(\mathbb{Z}[t],X)$-graph $\Gamma$ is called
{\em $U$-folded} if for any reduced path $p$ in $\Gamma$ with
$\overline{\mu(p)} = w$  there exists a unique label reduced path
$q$ such that $o(q) = o(p), t(q) = t(p), \mu(q) = \pi(w)$.

The above definition is equivalent to a more technical one given
in \cite{MRS2}.

\begin{prop} \cite{MRS2}
\label{pr:4.2.1} Let $\Gamma$ be a finite connected
$(\mathbb{Z}[t],X)$-graph and $v \in V(\Gamma)$. Then there exists a
$U$-folded $(\mathbb{Z}[t],X)$-graph $\Delta$ and $v' \in V(\Delta)$
such that $L(\Gamma,v) = L(\Delta,v')$. Moreover $\Delta$ can be constructed
effectively by adding to $\Gamma$ finitely many edges and applying finitely many
free and $U$-foldings.
\end{prop}

\begin{prop} \cite{MRS2}
\label{pr:4.3.1} Let $H$ be a finitely generated subgroup of
$F^{\mathbb{Z}[t]}$. Then there  exists a $U$-folded
$(\mathbb{Z}[t],X)$-graph $\Gamma$ and a vertex $v$ of $\Gamma$
such that $L(\Gamma,v) = H$.
\end{prop}

\begin{prop} \cite{MRS2}
\label{pr:4.3.2} There is an algorithm which, given finitely many
standard decompositions of  elements $h_1, \ldots, h_k$ from
$F^{\mathbb{Z}[t]}$, constructs a $U$-folded
$(\mathbb{Z}[t],X)$-graph $\Gamma$, such that $L(\Gamma,v) =
\langle h_1,\ldots,h_k \rangle$.
\end{prop}

The properties of $U$-folded graphs make it possible to solve the
membership problem in finitely generated subgroups of
$F^{\mathbb{Z}[t]}$.

\begin{prop} \cite{MRS2}
\label{pr:4.3.3} Every finitely generated subgroup of
$F^{\mathbb{Z}[t]}$ has a solvable membership  problem. That is,
there exists an algorithm which, given finitely many standard
decompositions of elements $g, h_1, \ldots, h_k$ from
$F^{\mathbb{Z}[t]}$, decides whether or not $g$ belongs to the
subgroup $H = \langle h_1,\dots , h_n\rangle$ of
$F^{\mathbb{Z}[t]}$.
\end{prop}

\section{Finite index criteria}
\label{sec:findex}

It is not difficult to check if a finitely generated subgroup $H$ of a free group $G$
is of finite index (see, for example, \cite{KM}). This can be done by
checking if the normal form of every element of $G$ is ``readable'' in
a folded graph $\Gamma_H$ corresponding to $H$. Similar result can be
easily proved for finitely generated subgroups of $F^{\mathbb{Z}[t]}$.

\begin{prop}
\label{fincrit1}
Let $G$ be a finitely generated subgroup of $F^{\mathbb{Z}[t]}$ and $H \leq G$. Then
the following are equivalent:
\begin{enumerate}
\item $|G : H| < \infty$,
\item there exists a finite $U$-folded ($\mathbb Z[t],X$)-graph $\Delta$ with a vertex
$v$ such that $H = L(\Delta,v)$ and for every $g \in G$ there exists a path $p$
in $\Delta$ such that $o(p) = v,\ \overline{\mu(p)} = g$.
\end{enumerate}
\end{prop}
{\it Proof.} At first, assume $|G : H| < \infty$. Hence, there exist $g_1,\ldots,g_k \in G$
such that
$$G = H \cup H g_1 \cup \cdots \cup H g_k.$$
Take a finite $U$-folded ($\mathbb Z[t],X$)-graph $\Gamma$ with a vertex $v$ such that
$H = L(\Gamma,v)$. For each $g_i,\ i \in [1,k]$ take a path labeled by $\pi(g_i)$ and
glue its initial end-point to $\Gamma$ at $v$. The resulting ($\mathbb Z[t],X$)-graph $\Gamma'$
by Proposition \ref{pr:4.2.1} can be transformed into a $U$-folded ($\mathbb Z[t],X$)-graph
$\Delta$ whose language is $H$. But since for every product $h \ast g_i,\ h \in H, i \in [1,k]$
there exists a path $p$ in $\Gamma'$ such that $o(p) = v,\ \overline{\mu(p)} = g$, this property also
holds in $\Delta$.

\smallskip

Now, assume that there exists a finite $U$-folded ($\mathbb Z[t],X$)-graph $\Delta$ with a vertex
$v$ such that $H = L(\Delta,v)$ and for every $g \in G$ there exists a path $p$
in $\Delta$ such that $o(p) = v,\ \overline{\mu(p)} = g$. For every $w \in \Delta$
there exists a path $p_w$ such that $o(p_w) = v,\ t(p_w) = w$. Since $\Delta$ is finite, the set of
such paths $p_w$ is finite and their reduced labels obviously can be taken to be representatives of
right cosets in $G$ by $H$.

\hfill $\square$

At the same time, it is important to understand that not every $U$-folded graph representing a
subgroup $H$ of finite index in $G$ has the property that the normal form of every element of $G$ is
``readable'' in it.

\begin{figure}[htbp]
\centering{\mbox{\epsfig{figure=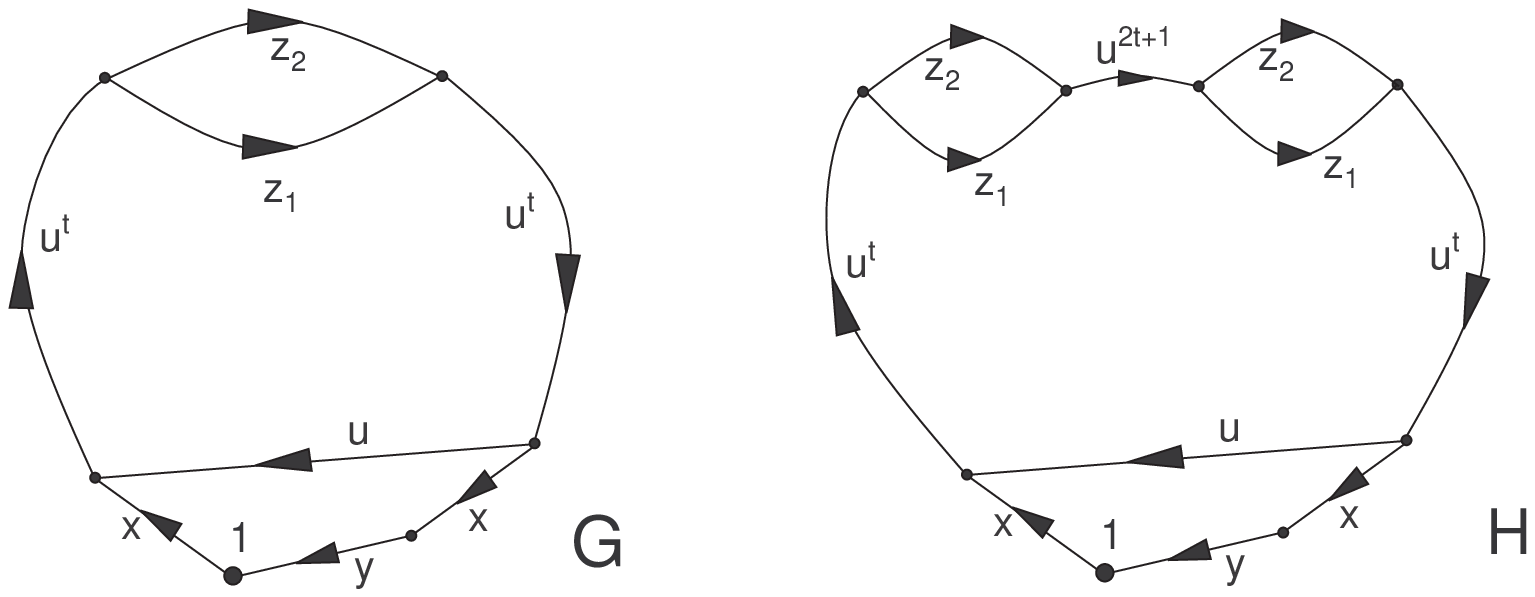,height=1.4in}}}
\caption{}
\label{graph1}
\end{figure}

\begin{example}
{\rm Let $u = x y x \in F(X) \cap U$ and $G = \langle a,b\rangle,\ H =
\langle a^2,b^2,ab\rangle$, where $a = x u^t z_1 u^t x y,\ b = x u^t z_2 u^t x y$,
Without loss of generality we can assume $z_1, z_2 \in F(X)$ to be such that
the graphs shown on Figure~\ref{graph1} are $U$-folded. Observe that $|G:H| = 2$,
but $g = x u^t z_1 u^t x y \in G$ is not ``readable'' in the graph defining $H$.}
\end{example}

Hence, ``readability'' of normal forms of elements from $G$ is a too strong property to work with.
Instead, below we develop the idea of ``readability'' of infinite paths arising in a $U$-folded graph
for $G$, and show that every such infinite path is readable in any $U$-folded graph for $H$ as long as
$|G : H| < \infty$.

\subsection{Equivalence of infinite powers}
\label{subs:w-equiv}

For any $K \subseteq \mathbb{Z}[t]$, denote
$$u^K = \{ u^\alpha \mid \alpha \in K\}.$$
Let $W$ be a finite subset of $U$. Denote $B_W = \{X \cup X^{-1}\} \cup (\bigcup_{u \in W} u^{\mathbb Z[t]})$.
For any $u \in W$ we say that {\it $u^\alpha$ is $W$-equivalent to $u^\beta$}, where $\alpha,\beta \in
\mathbb Z[t] - \mathbb Z$, and denote $u^\alpha \sim_W u^\beta$, if the following condition holds:
for any $w_1, w_2 \in B_W^*$ if $w_1 u^\alpha w_2$ is a standard form, then $w_1 u^\beta w_2$ is a
standard form as well, provided $w_1$ does not end with a power of $u$ and $w_2$ does not start with
a power of $u$.

\begin{prop}
\label{w-equiv}
For every finite $W \subset U$ and $u \in W,\ u^{\mathbb Z[t]}$ splits into a finite number of $W$-equivalence classes.
Moreover, for any $\alpha \in \mathbb{Z}[t]$, the equivalence class of $u^\alpha$ can be effectively constructed.
\end{prop}
{\it Proof.} Fix a finite set $W$ and $u \in W$.

For a given occurrence $w_1 u^\alpha w_2$, describe all $u^\beta$ that $u^\alpha$ can be replaced with,
not breaking the standard form. Suppose some $u^\beta$ does break the standard form. Enumerate possible reasons according
to definition of standard form:
\begin{enumerate}
\item $w_1 u^\beta w_2$ is no longer reduced, which means that $\beta$ is of opposite sign to $\alpha$,
\item $w_1 = t_1 s_1,\ w_2 = s_2 t_2,\ s_1 u^\beta s_2 = \pi(v^\kappa),\ v \in U, v > u$, and either
$t_1$ ends with $v^\delta$, or $t_2$ starts with $v^\delta$,
\item $w_1 = t_1 s_1,\ s_1 u^\gamma = \pi(v^\kappa)$, where either $\beta \ge \gamma > 0$ or $\beta \le \gamma < 0$, and either
$t_1$ ends with $v^\delta$, or $w_2$ starts with $v^\delta$,
\item $w_2 = s_2 t_2,\ u^\gamma s^2 = \pi(v^\kappa)$, where either $\beta \ge \gamma > 0$ or $\beta \le \gamma < 0$, and either
$t_1$ ends with $v^\delta$, or $w_2$ starts with $v^\delta$.
\end{enumerate}
As we can see, all possible cases result in conditions of the following types
$$\beta < \gamma,\ \beta \le \gamma,\ \beta > \gamma,\ \beta \ge \gamma,\ \beta = \gamma,\ \beta \neq \gamma.$$
Since there are only finitely many occurrences of $u$ in standard forms of $v \in W,\ v > u$, these
condition split $\mathbb Z[t] - \mathbb Z$ into finitely many classes and each class can be easily constructed.

\hfill $\square$

\begin{cor}
\label{cor:w-equiv}
For every finite $W \subset U$ and $u \in W$, $W$-equivalence classes in $u^{\mathbb Z[t]}$ can be effectively
described and enumerated.
\end{cor}

{\it Proof.} Follows immediately from the previous result.

\hfill $\square$

\subsection{Standard form types associated to ($\mathbb Z[t],X$)-graph}
\label{subs:types}

Let $\Gamma$ be a finite $U$-folded ($\mathbb Z[t],X$)-graph. We fix $\Gamma$ for the
rest of this subsection. Observe that $U(\Gamma)$ is a finite subset of $U$ and
for each $u \in U(\Gamma)$, by Proposition \ref{w-equiv}, there are only finitely many
$U(\Gamma)$-equivalence classes in $u^{\mathbb{Z}[t]}$.

Fix $u \in U(\Gamma)$. Let $C$ be a $u$-component of $\Gamma$. By Lemma \ref{le:3.3.1}, $H_u(C)$ is a subgroup of
$u^{\mathbb{Z}[t]}$. For every $a,b \in V(C)$ define the set $u(a,b)$ of reduced paths in $C$ from
$a$ to $b$ and a subset $H_u(a,b)$ of $u^{\mathbb{Z}[t]}$ as follows
$$H_u(a,b) = \{\overline{\mu(p)}\, \mid \, p \in u(a,b) \}.$$
Note that $H_u(a,b) = H_u(C) \ast \overline{\mu(p)}$ for any reduced path $p$ in
$C$ from $a$ to $b$. That is, $H_u(a,b)$ is a (right) coset in $u^{\mathbb{Z}[t]}$
by $H_u(C)$, and it is completely defined by a finite number of generators of $H_u(C)$
and a representative $\overline{\mu(p)}$.

For a $u$-component $C$ of $\Gamma$, a pair $a,b \in V(C)$, and a $U(\Gamma)$-equivalence class $u^A$ in
$u^{\mathbb{Z}[t]}$, consider a subset $u^A(a,b)$ of $u(a,b)$ defined as
$$u^A(a,b) = \{p \in u(a,b)\, \mid \, \overline{\mu(p)} \in u^A \}.$$
Paths $p$ and $q$ which belong to a certain $u^A(a,b)$ we call {\it $U(\Gamma)$-equivalent}.

\begin{prop}
\label{equiv-effective}
For every $u \in U(\Gamma)$, there are only finitely many sets of the type $u^A(a,b)$. Moreover, the set of
labels of paths from each $u^A(a,b)$ can be effectively described and enumerated.
\end{prop}
{\it Proof.} The first part of the statement follows from the fact that the number of tuples $(C,a,b,u^A)$ is finite
(in particular, by Proposition \ref{w-equiv}).

Next, by definition, the set of labels of paths from $u^A(a,b)$ is $H_u(a,b) \cap\, u^A$, where both sets can be
effectively described and enumerated (in particular, by Corollary \ref{cor:w-equiv}). Finally, since these sets are
subsets of an abelian group $u^{\mathbb{Z}[t]}$, the intersection also can be effectively described and enumerated.

\hfill $\square$

Denote by $\paths{\Gamma}$ the set of all paths in $\Gamma$ and define the set of {\it special paths}
as follows
$$\spaths{\Gamma} = \{ p \in \paths{\Gamma} \mid \mu(p) = \pi(g)\ {\rm for\ some}\ g \in
F^{\mathbb{Z}[t]}\}.$$
For a vertex $v \in V(\Gamma)$, similarly define
$$\spaths{\Gamma}_v = \{ p \in \spaths{\Gamma} \mid o(p) = v\}.$$
Observe that for every $p \in \spaths{\Gamma}$ there is a {\it standard decomposition}
$$p = p_1 \cdots p_k,$$
where either
\begin{enumerate}
\item $p_i$ is an edge in $\Gamma(0)$, or
\item $p_i \in u(a,b)$, for some $u$-component $C$ and $a,b \in V(C)$, is maximal
with respect to inclusion.
\end{enumerate}
From the above definition it is easy to draw the following result.

\begin{lemma}
\label{stand-unique}
Let $p_1 \cdots p_k$ and $q_1 \cdots q_n$ be standard decompositions of $p \in \spaths{\Gamma}$.
Then
\begin{enumerate}
\item $k = n$,
\item if $p_i \in E(\Gamma(0))$ then $q_i \in E(\Gamma(0))$ and $p_i = q_i$,
\item if $p_i \in u(a,b)$, for some $u$-component $C$ and $a,b \in V(C)$, then $q_i \in u(a,b)$
and $\overline{\mu(p_i)} = \overline{\mu(q_i)}$.
\end{enumerate}
\end{lemma}

Let $p,q \in \spaths{\Gamma}_v$ for some $v \in V(\Gamma)$ and let $p = p_1 \cdots p_k,\ q = q_1 \cdots q_n$
be some of their standard decompositions. We say $p$ and $q$ are {\it $U(\Gamma)$-equivalent} if
\begin{enumerate}
\item $k = n$,
\item $p_i = q_i$ if both are edges in $\Gamma(0)$,
\item $p_i, q_i \in u^A(a,b)$ for some $u$-component $C$, its vertices $a,b$, and a $U(\Gamma)$-equivalence
class $u^A$ in $u^{\mathbb{Z}[t]}$.
\end{enumerate}
Observe that $U(\Gamma)$-equivalent paths have the same end-points. $U(\Gamma)$-equivalence classes of
special paths we call {\it types}. The $U(\Gamma)$-equivalence class of a special path $p$ we denote $\type{p}$.

Let $p \in \spaths{\Gamma}$ and $p = p_1 \cdots p_n$ be a standard decomposition of $p$. Hence, $\type{p}$ can
be represented as $t_1 \cdots t_n$, where each $t_i$ is either an edge labeled by a letter from $X \cup X^{-1}$,
or $u^A(a,b)$ for some $u$-component $C$, its vertices $a,b$, and a $U(\Gamma)$-equivalence
class $u^A$ in $u^{\mathbb{Z}[t]}$. Uniqueness of such a representation follows from Lemma \ref{stand-unique} and
the definition of $U(\Gamma)$-equivalence of special paths. Observe that every type can be viewed as a finite word
in the alphabet $T_\Gamma = E(\Gamma(0)) \cup \{u^A(a,b)\}$, which is finite by Proposition \ref{equiv-effective}.
The length of $t$ as a finite word in the alphabet $T_\Gamma$ denote by $\|t\|$, that is $\|t\| = n$. Finally,
define the {\it label} of $t$ to be $\mu(t) = \mu(t_1) \cdots \mu(t_n)$, where
\[ \mu(t_i) = \left\{ \begin{array}{ll}
\mbox{$\mu(t_i)$,} & \mbox{if $t_i \in \Gamma(0)$} \\
\mbox{$u^A$,} & \mbox{if $t_i = u^A(a,b)$}
\end{array}
\right.
\]

Denote by $\types{\Gamma}\subset (T_\Gamma)^\ast$ the set of all types in $\Gamma$. Also,
denote by $\types{\Gamma}_v$ the set of all types in $\Gamma$ ``readable''
from $v \in V(\Gamma)$.

An infinite sequence of types $t_1, t_2, \ldots, t_n, \ldots$ in $\Gamma$, where $t_i$ is an
initial subword of $t_j$ for any $i < j$, is, naturally, called an {\em infinite type},
considering it as an element of $(T_\Gamma)^\omega$. As in the case of finite types we define
$\types{\Gamma}^\omega$ and $\types{\Gamma}_v^\omega,\ v \in V(\Gamma)$.

\begin{remark}
\label{subword}
$\types{\Gamma}$ is closed under taking subwords.
\end{remark}

\begin{lemma}
\label{decide_inf}
Given a $U$-folded ($\mathbb Z[t],X$)-graph $\Gamma$ it is possible to decide
effectively if $\types{\Gamma}$ is infinite.
\end{lemma}
{\it Proof.} If $\types{\Gamma}$ is finite then $\Gamma$ cannot contain more than
one nontrivial $u$-component for $u \in U$, and there cannot be any cycles in
$\Gamma(0)$ save for those doubled in the $u$-component. So, $L(\Gamma,v)$ is abelian
for any $v \in V(\Gamma)$. Observe that existence of $u$-components and loops in
$\Gamma(0)$ is algorithmically decidable.

\hfill $\square$

\begin{lemma}
\label{exist_inf}
$\ \types{\Gamma}$ is infinite if and only if $\types{\Gamma}^\omega$ is not empty.
\end{lemma}
{\it Proof.} Immediate from K\"{o}nig's Lemma, since $T_\Gamma$ is finite.

\hfill $\square$

Define a constant $M(\Gamma) \in \mathbb{N}$ as follows
$$M(\Gamma) = 1 + \max_{u\in U(\Gamma)} \|u^{\pm 1}\|.$$

\begin{lemma}
\label{tsr} For any $r,s,t \in \types{\Gamma}$ with $\|s\| \ge M(\Gamma)$ if
$ts,sr \in \types{\Gamma}$ then $t s r \in \types{\Gamma}$.
\end{lemma}
{\it Proof.} Suppose $r,s,t \in \types{\Gamma}$ and let $p_r, p_s, p_t \in \spaths{\Gamma}$
be such that $r = \type{p_r},\ s = \type{p_s},\ t = \type{p_t}$. Since $t s, s r \in
\types{\Gamma}$ it follows that $t s = \type{p_t p_s}$ and $s r = \type{p_s p_r}$.
We are going to show that if $s$ is long enough then $\mu(p_t p_s p_r) =
\pi(\overline{\mu(p_t p_s p_r)})$.

Suppose this is not the case. Let
$$p_t p_s p_r = p_1 \cdots p_n,$$
where either $p_i \in \Gamma(0)$, or $p_i$ is a maximal (with respect to inclusion)
subpath inside of a $u$-component of $\Gamma$. From the definition of standard form it
follows that there is $p_k$ inside of a $u$-component of $\Gamma$ and a subpath $p'$ of
$p_t p_s p_r$ such that either
\begin{enumerate}
\item[(a)] $p' = p_i \cdots p_{k-1},\ \mu(p') = \pi(u^{\pm 1})$, or
\item[(b)] $p' = p_i \cdots p_{k-1},\ \overline{\mu(p')} = w \circ u^{\pm 1}$ and $\mu(p_i) = v^\alpha, v < u$, or
\item[(c)] $p' = p_{k+1} \cdots p_i,\ \mu(p') = \pi(u^{\pm 1})$, or
\item[(d)] $p' = p_{k+1} \cdots p_i,\ \overline{\mu(p')} = u^{\pm 1} \circ w$ and $\mu(p_i) = v^\alpha, v < u$.
\end{enumerate}
Assume $\|s\| > M(\Gamma)$, that is, $\|s\|>\|u\|$ for any $u \in U(\Gamma)$. It follows that
\begin{enumerate}
\item if $p_k \in p_t$ then $p'$ is a subpath of $p_t p_s$,
\item if $p_k \in p_s$ then $p'$ is a subpath of either $p_t p_s$, or $p_s p_r$,
\item if $p_k \in p_r$ then $p'$ is a subpath of $p_s p_r$.
\end{enumerate}
Since $t s$ and $s r$ are types it follows that we get a contradiction in any of the above cases.

\hfill $\square$

\begin{lemma}
\label{discard}
If $t_1, t_2, t_3, s \in \types{\Gamma}$ and $\|s\| \ge M(\Gamma)$ then $t_1 s t_2 s t_3 \in \types{\Gamma}$
if and only if $t_1 s t_3,\ s t_2 s \in \types{\Gamma}$.
\end{lemma}
{\it Proof.} Follows from Remark \ref{subword} and Lemma \ref{tsr}.

\hfill $\square$

Lemma \ref{discard} explains how to discard a loop from a standard form.

\begin{lemma}
\label{repeat}
Let $g \in F^{\mathbb{Z}[t]}$ be such that $\pi(g^n)= \pi(g)^n$ for any $n > 0$
and let $p$ be a path in $\Gamma$ such that $\mu(p) = \pi(g)$. Then for any
$h \in F^{\mathbb{Z}[t]}$
\begin{enumerate}
\item either there exists $k \in \mathbb{N}$ such that for all $n > k$,
$\pi(g^n \ast h) = \pi(g^{n-k}) \pi(g^k \ast h)$, or
\item $g = g_1 \circ g_2,\ p = p_1 p_2$ with $\mu(p_1) = \pi(g_1), \mu(p_2) = \pi(g_2)$,
where $g_2 \circ g \circ g_1 = u^\alpha$ for some $u \in U$, and $h = g_1 \circ u^\beta
\circ h'$,
\end{enumerate}
and, similarly,
\begin{enumerate}
\item either there exists $k \in \mathbb{N}$ such that for all $n > k$,
$\pi(h \ast g^n) = \pi(h \ast g^k) \pi(g^{n-k})$, or
\item $g = g_1 \circ g_2,\ p = p_1 p_2$ with $\mu(p_1) = \pi(g_1), \mu(p_2) = \pi(g_2)$,
where $g_2 \circ g \circ g_1 = u^\alpha$ for some $u \in U$, and $h = h' \circ u^\beta
\circ g_2$.
\end{enumerate}
\end{lemma}
{\it Proof.} The word part of statement is established
by inspecting heights of $h,g$. Then the path part of statement
holds, since division $g = g_1 \circ g_2$ cannot correspond to
subdividing a $u$-edge of $p$ (if that were the case, $\pi(u^{2\alpha})$ would be not
equal to $\pi(u^\alpha)^2$).

\hfill $\square$

An infinite type $t \in \types{\Gamma}^\omega$ is called {\em (infinite) periodic}
if $t = t_1 t_2 t_2 \cdots = t_1 t_2^\infty$, where $t_1,\ t_2,\ t_1 t_2 \in
\types{\Gamma}$.

Define
$$N(\Gamma) = M(\Gamma)(1 + |T_\Gamma|^{M(\Gamma)}).$$

\begin{prop}
\label{exist_period}
If $\ \types{\Gamma}^\omega$ is non-empty then there exists a periodic type
$t = t_1 t_2^\infty$ in $\Gamma$. Moreover, $\|t_1\| + \|t_2\| < N(\Gamma)$.
\end{prop}
{\it Proof.} Consider $t = \tau_1 \tau_2 \cdots \tau_k \cdots \in \types{\Gamma}^\omega$.
Hence, $\tau_1 \cdots \tau_{N(\Gamma)}$ contains at least two non-intersecting copies of a
subword $s$ of length at least $M(\Gamma)$, that is,
$$t = t_1 s t_2 s t',\ \|t_1 s t_2 s\| \le N(\Gamma),\ t' \in \types{\Gamma}^\omega.$$
Then by Lemma \ref{tsr} applied to $(t_1 s t_2) s$ and $s (t_2 s)$ we have
$t_1 (s t_2)^\infty \in \types{\Gamma}^\omega$ is periodic.

\hfill $\square$

We say that a periodic type $t_1 t_2^\infty$ is of {\em content} $K \in \mathbb{N}$
if $\|t_1\| + \|t_2\| \le K$.

\begin{remark}
\label{enumerable}
For any fixed $K\in \mathbb{N}$ one can effectively enumerate all periodic types of
content $K$.
\end{remark}

\subsection{Finite index conditions}
\label{subs:fi}

Let $G$ and $H$ be finitely generated subgroups of $F^{\mathbb{Z}[t]}$
such that $H \le G$. Let $G = L(\Gamma, 1_G),\ H = L(\Delta, 1_H)$, where $\Gamma$ and $\Delta$
are $U$-folded ($\mathbb Z[t],X$)-graphs, $1_G \in V(\Gamma), 1_H \in V(\Delta)$.
In this subsection we describe an algorithm which decides if $|G : H| <\infty$.

\medskip

Without loss of generality we can assume $\types{\Gamma}^\omega$ to be non-empty.
Observe that otherwise $G$ is abelian and it is easy to check if $|G : H| <\infty$.

\medskip

A path $p \in \paths{\Gamma}$ is \emph{doubled} in $\Delta$ at $v \in V(\Delta)$
if there is a path $\hat{p} \in \paths{\Delta}_v$ such that $\mu(\hat{p}) =
\pi(\overline{\mu(p)})$.

A type $t \in \types{\Gamma} \cup \types{\Gamma}^\omega$ is
\begin{itemize}
\item \emph{doubled} in $\Delta$ at $v \in V(\Delta)$ if there exists $\hat{t} \in \types{\Delta}_v$
such that $\mu(t) \subseteq \mu(\hat{t})$,
\item \emph{$N$-almost doubled} in $\Delta$ at $v \in V(\Delta)$ if $t = t_1 t_2$ is a finite type
for which $t_1$ is doubled in $\Delta$ at $v \in V(\Delta)$ and $\|t_2\| \le N$.
\end{itemize}
Usually, we use $\hat{\ }$-symbol to denote a doubling type (path). Also, below we do not specify
a point at which the type is doubled when it is clear from the context.
For example, if a type $t = t_1 t_2 \in \types{\Gamma} \cup \types{\Gamma}^\omega$ is doubled
at $\Delta$ at $v \in V(\Delta)$ then $t_2$ is doubled in $\Delta$ at $v' \in V(\Delta)$,
where $v' = t(p),\ p \in \spaths{\Delta},\ \mu(p) \in \mu(t_1)$.

Define
$$K(\Gamma,\Delta) = M(\Gamma)(|T_\Gamma|^{M(\Gamma)}(|V(\Delta)|-1) + 1).$$

\begin{lemma}
\label{long}
If all periodic types in $\types{\Gamma}^\omega_v$ of content $K = K(\Gamma,\Delta)$ are doubled
in $\Delta$ at $w \in V(\Delta)$ then any type $t \in \types{\Gamma}_v$
of length greater than $K$ (or $t \in \types{\Gamma}^\omega_v$) can be represented as follows
$$t = t^{(1)} s^{(1)} t^{(2)} s^{(2)} t^{(3)},$$
where
\begin{itemize}
\item $s^{(1)} = s^{(2)} = s,\ \|s\| \ge M(\Gamma)$,
\item $t^{(1)}, s^{(1)}, t^{(2)}, s^{(2)}$ are doubled in $\Delta$ by
$\hat t^{(1)},\hat s^{(1)},\hat t^{(2)},\hat s^{(2)}$ respectively,
\item $\hat s^{(1)}=\hat s^{(2)}=\hat s$.
\end{itemize}
\end{lemma}
{\it Proof.} Assume that all periodic types in $\types{\Gamma}^\omega_v$ of content $K$ are doubled in $\Delta$
at $w \in V(\Delta)$.

Let $t \in \types{\Gamma}_v$ be of length greater than $K$ (or $t \in \types{\Gamma}^\omega_v$).
Observe that $t$ contains at least $|V(\Delta)|$ non-intersecting copies of a subword of length at least $M = M(\Gamma)$.
That is, $t = t_1 s_1 t_2 s_2 \cdots t_k s_k t_{k+1},\ s_i = s$ with $\|t_1 \cdots s_k\| \le K$.
Note that by Lemma \ref{tsr}, $t_\infty = t_1 (s_1\cdots t_k)^\infty$ is a periodic type
of content at most $K$. By our assumption, $t_\infty$ is doubled in $\Delta$ by a type
$\hat t_\infty = \hat t_1 \hat s_1 \hat t_2 \hat s_2 \cdots \hat t_k \hat s_k \cdots$,
where $\hat t_i$ and $\hat s_i$ are doubles of $t_i$ and $s_i$, respectively. Hence, we have
$\hat s_i = \hat s_j$ for some $i,j$, so we set $t^{(1)} = t_1 s_1 \cdots t_i,\
s^{(1)} = s_i,\ t^{(2)} = t_{i+1} s_{i+1} \cdots t_j,\ s^{(2)}=s_j,\ t^{(3)} = t_{j+1} s_{j+1} \cdots t_{k+1}$.

\hfill $\square$

\begin{cor}
\label{co:long}
If all periodic types in $\types{\Gamma}^\omega_v$ of content $K = K(\Gamma,\Delta)$ are doubled
in $\Delta$ at $w \in V(\Delta)$ then
\begin{enumerate}
\item any $t \in \types{\Gamma}_v,\ \|t\| \ge K$ can be represented as $t = t' t'' t''',\ \|t''\| > 0$,
where $t' t''' \in \types{\Gamma}_v$, and $t$ is $K$-almost doubled if and only if $t' t'''$ is
$K$-almost doubled,
\item  any $t \in \types{\Gamma}^\omega_v$ can be represented as $t = t' t'' t''',\ \|t''\| > 0$,
where $t' t''' \in \types{\Gamma}^\omega_v$, and $t$ is doubled if and only if $t' t'''$ is
doubled.
\end{enumerate}
\end{cor}
{\it Proof.} Follows from Lemma \ref{long} by setting $t' = t^{(1)} s^{(1)},\ t'' =
t^{(2)} s^{(2)},\ t''' = t^{(3)}$.

\hfill $\square$

If $t = t' t'' t'''$ is the decomposition from the corollary above then we call
$t' t'''$ a {\em reduction} of $t$. Respectively, $t$ is a {\em lift} of $t' t'''$.

\begin{lemma}
\label{periodic}
If all periodic types in $\types{\Gamma}^\omega_v$ of content $K = K(\Gamma,\Delta)$ are doubled
in $\Delta$ at $w \in V(\Delta)$ then all types in $\types{\Gamma}_v$ are $K$-almost doubled
in $\Delta$ at $w \in V(\Delta)$. Moreover, all types in $\types{\Gamma}^\omega_v$ are doubled.
\end{lemma}
{\it Proof.} Let $t \in \types{\Gamma}_v$. Consider a sequence
$$t_1, t_2, \ldots, t_{k+1}, \ldots$$
where $t_1 = t,\ t_{i+1} = t'_i t'''_i$ is a reduction of $t_i$, and $\|t_{k+1}\|< K$ (such $k$
exists since $\|t\|$ is finite). Hence, $t_{k+1}$ is $K$-almost doubled implying that $t_k$ is
$K$-almost doubled, and after $k$ lifts we get $t$ which is $K$-almost doubled.

If $t \in \types{\Gamma}^\omega_v$ then $t$ can be viewed as a sequence $\{t_i\}$ of finite types of
increasing length which extend each other. By the argument above, $t_i = s_i r_i$, where $s_i$ is
doubled and $\|r_i\| \le K$. Hence, $\|s_i\|$ increases with growth of $i$ implying that all
$t_i$ are doubled. So, $t$ is doubled as well.

\hfill $\square$

In the next proposition we elaborate on the following idea. Suppose
there is a type $t \in \types{\Gamma}_v$ which is not doubled in $\Delta$.
Whether we can use it to produce infinitely many cosets in $G$ by $H$ depends
on whether or not we can extend $t$ into an infinite type. The former implies the
index is infinite. If the latter holds for all non-doubled types then all possible non-doubled
pieces occur in the ``ends'' of these types (otherwise we can extend them without
changing non-doubled pieces!). Lemma \ref{periodic} puts a bound on how long the ``end''
can be, virtually explaining why we can cover all $G$ with a finite number of cosets by $H$.

\begin{prop}
\label{index}
The following statements are equivalent:
\begin{enumerate}
\item $|G:H| < \infty$,
\item each periodic type in $\types{\Gamma}^\omega_{1_\Gamma}$ of content $K(\Gamma,\Delta)$
is doubled in $\Delta$ at $1_H$,
\item each infinite type in $\types{\Gamma}^\omega_{1_\Gamma}$ is doubled in $\Delta$ at $1_H$.
\end{enumerate}
\end{prop}
{\it Proof.} (2) $\Leftrightarrow$ (3) follows from Lemma \ref{periodic}.

\smallskip

(1) $\Rightarrow$ (2) Suppose there exists $t = t_1 t_2^\infty \in \types{\Gamma}^\omega_{1_\Gamma}$
which is not doubled in $\Delta$ at $1_H$. Let $q_1, q_2 \in \spaths{\Gamma}$ such that
$t_1 = \type{q_1},\ t_2 = \type{q_2}$, and $q_1q_2^\infty$ is not doubled in $\Delta$.
Since by the assumption $|G:H| < \infty$,
there exists a finite set $g_1, g_2, \ldots, g_m$ in $G$ such that for any $g \in G$ there exists
$i \in [1,m]$ such that $g \ast g_i \in H$. We show that there exists a loop in
$\Gamma$ at $1_G$ corresponding to some $g \in G$ such that none of standard
forms $\pi(g \ast g_i)$ is readable in $\Delta$ from $1_H$.

Fix $i \in [1,m]$. Let $\pi(g_i) = \mu(s_i),\ s_i \in \spaths{\Gamma}_{1_G}$ and paths
$r_n = q_2^n q_1^{-1} s_i$. By Lemma \ref{repeat} there are two possibilities.
\begin{enumerate}
\item[(a)]
There exists $k(i)$ such that for any $n > k(i),\ \type{r_n}$ begins with $t_2^{n-k(i)}$.
Then by Lemma \ref{tsr}, $t_1 \cdot \type{r_n} \in \types{\Gamma}_{1_\Gamma}$ and its initial
subtype $t_1 t_2^{n-k(i)}$ is not doubled. Hence, for any path $r \in \spaths{\Gamma}$ such
that $\type{r} = \type{r_n}$ we can choose a path $p$ such that $\mu(p) \in \mu(t_1t_2^n)$
such that $\mu(p) \mu(r)$ is a standard form of an element of $G$, and $\mu(p)$ cannot be read in
$\Delta$. This implies that $\overline{\mu(p) \mu(r)} \notin H$.

\item[(b)]
$q_2 = q' q'',\ \overline{\mu(q'' q_2^l q')} = u^\kappa,\ u \in U(\Gamma)$ and
$\overline{\mu(q_1^{-1} s_i)} = \overline{\mu(q')} \circ u^\beta \circ w$. Let
$s = \type{q_1 q_2^n q_1^{-1} s_i}$ $= t' u_a t''$, where
$u_a\in T_\Gamma$ corresponds to the designated
occurrence of $u^\beta$, and $o(u_a)=a$.

Suppose $s$ is doubled in $\Delta$ at $1_H$. Then $H_u(a) \cap \langle u\rangle\neq\varepsilon$ since
$t_1 t_2^n \in \types{\Gamma}_{1_\Gamma}$ for any $n$, and $H_u(a')\cap \langle
u\rangle = \varepsilon$ ($a'$ corresponds to $a$ in the double of $s$) since $q_1 q_2^\infty$ is not
doubled. It follows that we can read only bounded finite powers of $u$ in $\comp(a')$ at $a'$, and
we can choose large enough $k(i)$ so that $q_1 q_2^n q_1^{-1} s_i, n>k(i),$ is not doubled in $\Delta$ at $1_H$.
\end{enumerate}

Taking $k = \max_{i=1}^m \{k(i)\}$ we obtain that $q_1 q_2^n q_1^{-1} s_i$ represents an element of $G$,
and is not doubled in $\Delta$ at $1_H$ for any $n > k$ and $i \in [1,m]$.

\smallskip

(2) $\Rightarrow$ (1) Suppose each periodic type of content at most $K(\Gamma,\Delta)$ is doubled.
We prove that $|G:H| < \infty$.

It is enough to construct a finite set of paths $r_1,\ldots, r_m \in \paths{\Gamma}_{1_G}$ such that
for any loop $p \in \spaths{\Gamma}_{1_G}$ such that $t(p) = 1_\Gamma$ there exists $i \in [1,m]$ for
which $\type{p r_i}$ is doubled in $\Delta$ at $1_H$. In this case the number of cosets in $G$ by
$H$ is at most $|V(\Delta)|$.

Observe that $t = \type{p}$ is $K(\Gamma,\Delta)$-almost doubled by Lemma \ref{periodic}.
Moreover, from the proof of Lemma \ref{periodic} we have $t = t_1 t_2$, where $t_1$ is a lift of $t' \in
\types{\Gamma}_{1_G}$ (both $t_1$ and $t'$ are doubled in $\Delta$ at $1_H$) and $\|t' t_2\| < K(\Gamma,\Delta)$.
Let $t_1 = \type{p_1},\ t_2 = \type{p_2},\ t' = \type{p'}$, where $p_1,\ p_2,\ p' \in \spaths{\Gamma}$. Observe that $t(p_1) = t(p')$.
It follows that $\type{p_1 p'^{-1}}$ is doubled in $\Delta$ at $1_H$. Indeed, $p'$ is doubled in $\Delta$ at $1_H$, so
$\type{p'^{-1}}$ is also doubled at $v \in V(\Delta)$ corresponding to the end of the double of $p'$, since both $\Gamma$ and
$\Delta$ are $U$-folded.

It is only left to note that $\overline{\mu(p_1 p'^{-1})} = \overline{\mu(p(p'p_2)^{-1})}$, and $p'p_2$ is a loop, so
$\mu((p'p_2)^{-1}) \subset G$. Since the number of types in $\Gamma$ of length at most $K(\Gamma,\Delta)$ is finite
the required statement follows.

\hfill $\square$

\begin{lemma}
\label{subgraph}
Any $t \in \types{\Delta}^\omega_{1_H}$ is doubled in $\Gamma$ at $1_G$.
\end{lemma}
{\it Proof.} Immediately follows since $H \le G$.

\hfill $\square$

\begin{lemma}
\label{check}
Let $t = t_1 t_2^\infty \in \types{\Gamma}^\omega_{1_G}$. Then there is an algorithm
which decides if $t$ is doubled in $\Delta$ at $1_H$.
\end{lemma}
{\it Proof.} We ``read'' $\mu(t_1 t_2^\infty)$ in $\Delta$ letter by letter. Observe that each letter
either is a letter from $X \cup X^{-1}$, or a set of the type $H_u(a,b) \cap\, u^A,\ u \in U(\Gamma)$,
which is effectively described. Hence, it is decidable if a letter
of $t$ can be ``read'' at a vertex of $\Delta$. Since $\Delta$ is a finite $U$-folded graph then either this ``reading''
fails at some vertex and $t$ is not doubled in this case, or some vertex of $\Delta$ is hit twice, and the doubling of
$t$ in $\Delta$ is found.

\hfill $\square$

\begin{theorem}
\label{index_main}
There is an algorithm which effectively decides if $|G:H| < \infty$.
\end{theorem}
{\it Proof.} By Lemma \ref{decide_inf} we can effectively determine
if $\types{\Gamma}, \types{\Delta}$ are finite.

If $\types{\Gamma}$ is finite then both $G$ and $H$
are abelian, so it is decidable if $|G:H| < \infty$.

If $\types{\Gamma}$ is infinite but $\types{\Delta}$ is finite then $|G:H| = \infty$
because $G$ contains a non-abelian free group and $H$ is abelian.

If $\types{\Gamma}$ and $\types{\Delta}$ are infinite then it is possible to check
the condition (2) of Proposition \ref{index} using Remark \ref{enumerable} and Lemma
\ref{check}.

\hfill $\square$

\begin{cor}
\label{cor:main1}
Let $G$ be a finitely generated subgroup of $F^{\mathbb{Z}[t]}$ and let $H \le G$ be finitely
generated. Then there is an algorithm which effectively decides if $|G:H|$ is finite.
\end{cor}

\begin{cor}
\label{cor:main2}
Let $G$ be a finitely generated fully residually free group and let $H \le G$ be finitely
generated. Then there is an algorithm which effectively decides if $|G:H|$ is finite.
\end{cor}

\section{Greenberg-Stallings Theorem}
\label{sec:green-stall}

In this section we prove the analog of the result for free groups
known as Greenberg-Stallings Theorem (see, for example, \cite{KM} Corollary 8.8).

Let $G_1, G_2$ be finitely generated non-abelian subgroups of $F^{\mathbb{Z}[t]}$ such that
$H = G_1 \cap G_2$ is of finite index in both $G_1$ and $G_2$. Let $\Gamma_1,\
\Gamma_2$, and $\Delta$ be $U$-folded graphs such that
$$G_1 = L(\Gamma_1, 1_{\Gamma_1}),\ G_2 = L(\Gamma_2, 1_{\Gamma_2}),\ H = L(\Delta, 1_\Delta),$$
where $1_{\Gamma_i} \in V(\Gamma_i),\ i = 1,2,\ 1_\Delta \in V(\Delta)$.
Observe that by Proposition \ref{index} every infinite type readable in $\Gamma_i$ at $1_{\Gamma_i}$,
for $i = 1,2$, is doubled in $\Delta$ at $1_\Delta$. At the same time, since
$$L(\Delta, 1_\Delta) = L(\Gamma_1, 1_{\Gamma_1}) \cap L(\Gamma_2, 1_{\Gamma_2})$$
it follows that every infinite type readable in $\Delta$ at $1_\Delta$ is doubled in
$\Gamma_i$ at $1_{\Gamma_i}$ for $i = 1,2$. Hence, it follows that every infinite type
readable in $\Gamma_1$ at $1_{\Gamma_1}$ is doubled in $\Gamma_2$ at $1_{\Gamma_2}$ and vice
versa. In particular, it follows that $U(\Gamma_1) = U(\Gamma_2)$.

\begin{lemma}
\label{alt-prod}
For any $g \in \langle G_1, G_2 \rangle$ the intersection $G_j \cap H^g$ is of finite index both in
$G_j,\ j = 1,2$ and $H^g$.
\end{lemma}
{\it Proof.} Any $g \in \langle G_1, G_2 \rangle$ can be viewed as a product $g = g_1 \cdots g_n$,
where $g_i, g_{i+1},\ i \in [1,n-1]$ do not belong to the same $G_j,\ j = 1,2$,
We use the induction on the lentgh of this product.

Let $g = f g_n$, where $f$ is a product of length $n-1$. By induction hypothesis $G_j \cap H^f$ is of
finite index both in $G_j,\ j = 1,2$ and $H^f$.

Suppose $g_n \in G_1$. Then obviously $(G_1)^{g_n} \cap (H^f)^{g_n} = G_1 \cap H^g = H^g$ is of finite index in $G_1$.
Next, $G_1 \cap G_2$ is of finite index in $G_1$, thus, $(G_1 \cap G_2) \cap H^g = G_2 \cap H^g < G_1$ is of finite
index in $G_1$. It follows that $G_2 \cap H^g < G_2$ is of finite index in $H = G_1 \cap G_2$, and, hence, in $G_2$

The same argument can be applied when $g = g_1 \in G_2$.

\hfill $\square$

\begin{cor}
\label{non-ab}
For any $g \in \langle G_1, G_2 \rangle$ the intersection $G_j \cap H^g$ is not abelian.
\end{cor}
{\it Proof.} Let $K < G < F^{\mathbb{Z}[t]}$. If $K$ is abelian and $|G : K| < \infty$ then $G$ is abelian too.
By our assumption both $G_1$ and $G_2$ are non-abelian, hence $H^g$ is non-abelian too for any
$g \in \langle G_1, G_2 \rangle$.

\hfill $\square$

Further we need the following corollary of Proposition 7 \cite{KMRS}.

\begin{lemma}
\label{double_cosets}
Let $G$ be a finite generated subgroup of $F^{\mathbb{Z}[t]}$ and $H,K \le G$ be finitely generated. Then
there exists a finite list of double cosets $K g_1 H, \ldots, K g_n H$ in $G$, such that for any $g \in G$ if
$g H g^{-1} \cap K$ is non-trivial and non-abelian then there exists $i \in [1,n]$ such that $g \in K g_i H$.
Moreover, $g_1,\ldots, g_n$ can be found effectively.
\end{lemma}
{\it Proof.} According to \cite{KMRS} (Proposition 7, cases 1 and 3), there exists a finite list of double cosets
$K g_1 H, \ldots, K g_n H$ in $F^{\mathbb{Z}[t]}$, such that for any $g \in G$ if $g H g^{-1}\cap K$ is
non-trivial and non-abelian then there exists $i \in [1,n]$ such that $g \in K g_i H$. Moreover, $g_1,\ldots, g_n$
can be found effectively.

Observe that we can effectively determine if $G \cap K g_i = \emptyset$ for every $i \in [1,n]$, and
$G \cap K g_i = \emptyset$ if and only if $G \cap K g_i H = \emptyset$. Hence, we can assume that for every
$g_i$ there exists at least one element $g \in G$ such that $K \cap H^g \neq 1$ and $g = f g_i h$ for some
$f \in K,\ h \in H$. It makes it possible to assume $g_1, \ldots, g_n \in G$.

\hfill $\square$

By the above lemma, there exists a finite list of double cosets $G_1 g_1 H, \ldots, G_1 g_n H$ in $\langle G_1, G_2 \rangle$,
such that for any $g \in \langle G_1, G_2 \rangle$ (since $G_1 \cap H^g$ is non-abelian for any such $g$) there exists $i \in [1,n]$
such that $g \in G_1 g_i H$. Since $|G_1 : H| < \infty$ we can use the list of double cosets
$H f_1 H, \ldots, H f_k H$ in $\langle G_1, G_2 \rangle$ by $H$.

Consider a $U$-folded graph $\Phi$ obtained from $\Delta$ and paths $p_1,\ldots,p_k$ labeled by $\pi(f_1),\ldots,\pi(f_k)$
whose initial end-points are glued to $1_\Delta$. Denote $v_i = t(p_i),\ i \in [1,k]$. The image of $1_\Delta$ in
$\Phi$ we still denote by $1_\Delta$. Obviously, $H = L(\Phi,1_\Delta)$.

Next, by our assumption for each $f_i$ there exists an element $g \in \langle G_1, G_2 \rangle$ such that
$g = h_1 f_i h_2$ for some $h_1,h_2 \in H$. Since $G_1 \cap H^g$ is of finite index in $G_1$ then
$G_1 \cap H^{f_i}$ is of finite index in $G_1$. Now, since both $H$ and $H^{f_i}$ have finite index in
$G_1$ it follows that $H \cap H^{f_i}$ has finite index in $G_1$ and hence in $H$. Let $w^{(i)}_1,\ldots,w^{(i)}_{m(i)}$
be a finite set of left coset representatives in $H$ by $H^{f_i}$.

Let $\Psi$ be a $U$-folded graph obtained from $\Phi$ and paths $p^{(1)}_1,\ldots,p^{(i)}_{m(i)}$ labeled by
$\pi(w^{(i)}_1),\ldots,$ $\pi(w^{(i)}_{m(i)})$, whose initial end-points are glued to $v_i$ for each $i \in [1,k]$.
The image of $1_\Delta$ in $\Psi$ we still denote by $1_\Delta$. Obviously, $H = L(\Psi,1_\Delta)$.

\begin{lemma}
\label{psi}
For any $g \in \langle G_1, G_2 \rangle$ there exists a path $p$ in $\Psi$ such that $o(p) = 1_\Delta,\
\mu(p) = \pi(g)$.
\end{lemma}
{\it Proof.} Take $g \in \langle G_1, G_2 \rangle$. Hence, it can be represented as $g = h_1 f_i h_2$ for some $f_i$ and
$h_1,h_2 \in H$. It follows that it can be represented further as $g = h_1 f_i w w^{(i)}_j$, where $w \in H^{f_i}$ and
$w^{(i)}_j$ is a coset representative in $H$ by $H^{f_i}$. By our construction there is a loop at $1_\Delta$ corresponding
to $h_1$, a path from $1_\Delta$ to $v_i$ corresponding to $f_i$, a loop at $v_i$ corresponding to $w$, and a path $p^{(1)}_j$
corresponding to $w^{(1)}_j$. Since $\Psi$ is folded there is also a required path $p$.

\hfill $\square$

\begin{cor}
$|\langle G_1, G_2 \rangle : H| < \infty$.
\end{cor}
{\it Proof.} Follows from Lemma \ref{psi} and Proposition \ref{fincrit1}.

\hfill $\square$

\begin{theorem}
\label{th:3.1}
Let $G_1, G_2$ be finitely generated subgroups of $F^{\mathbb{Z}[t]}$. If $H \leq G_1 \cap G_2$
is finitely generated and $|G_1 : H| < \infty,\ |G_2 : H| < \infty$ then $|\langle G_1, G_2 \rangle : H|
< \infty$.
\end{theorem}

\begin{theorem}
\label{th:3.2}
Let $G$ be a finitely generated fully residually free group and $G_1, G_2$ be finitely generated subgroups of $G$.
If $H \leq G_1 \cap G_2$ is finitely generated and $|G_1 : H| < \infty,\ |G_2 : H| < \infty$ then
$|\langle G_1, G_2 \rangle : H| < \infty$.
\end{theorem}

\section{Commensurator}
\label{sec:commensurator}

Let $G$ be a group and let $H \le G$. The {\em commensurator}
$\comm{G}{H}$ of $H$ in $G$ is defined as
$$\comm{G}{H} = \{g \in G \mid |H:H \cap g H g^{-1}| < \infty\ \mathrm{and\
} |gHg^{-1}:H\cap gHg^{-1}|<\infty \}.$$

It is easy to see that $\comm{G}{H}$ is a subgroup of $G$
containing $H$.

\begin{lemma}
\label{commensurator}
Let $G$ and $H$ be finitely generated non-abelian subgroups of $F^{\mathbb{Z}[t]}$
such that $H \le G$. Then there exist $f_1, \ldots, f_k \in G$ such that
$\comm{G}{H} = \langle H, g_1, \ldots, g_k \rangle$. Moreover, the elements $g_1, \ldots, g_k$
can be found effectively.
\end{lemma}
{\it Proof.} By Lemma \ref{double_cosets}, there exists a finite list of double cosets
$H g_1 H, \ldots, H g_n H$ in $G$, such that for any $g \in G$ if $g H g^{-1}\cap H$ is
non-trivial and non-abelian then there exists $i \in [1,n]$ such that $g \in H g_i H$. Moreover, $g_1,\ldots, g_n$
can be found effectively.

By Corollary \ref{cor:main1}, for each $g_i$ one can check effectively if
$|H : H \cap g_i H g_i^{-1}| < \infty$ and $|g_i H g_i^{-1} : H \cap g_i H g_i^{-1}| < \infty$,
so, we take only those elements $g_i$ for which both indeces are finite. Such elements form the required list
$f_1, \ldots, f_k$.

\hfill $\square$

\begin{theorem}
\label{effcomm}
Let $G$ and $H$ be finitely generated non-abelian subgroups of $F^{\mathbb{Z}[t]}$
such that $H \le G$. Then $\comm{G}{H}$ is finitely generated, and its generating set
can be found effectively.
\end{theorem}
{\it Proof.} Follows from Lemma \ref{commensurator}.

\hfill $\square$

\begin{cor}
\label{cor:effcomm}
Let $G$ be a finitely generated fully residually free group and let $H \le G$ be finitely
generated. Then $\comm{G}{H}$ is finitely generated, and its generating set
can be found effectively.
\end{cor}

\begin{lemma}
\label{normalizer}
Let $G$ and $H$ be finitely generated non-abelian subgroups of $F^{\mathbb{Z}[t]}$
such that $H \le G$. Then $|N_G(H) : H| < \infty$.
\end{lemma}
{\it Proof.} By Lemma \ref{double_cosets}, there exists a finite list of double cosets
$H g_1 H, \ldots, H g_n H$ in $G$, such that for any $g \in G$ if $g H g^{-1}\cap H$ is
non-trivial and non-abelian then there exists $i \in [1,n]$ such that $g \in H g_i H$. Moreover, $g_1,\ldots, g_n$
can be found effectively. So, if $S \subset \{g_1,\ldots, g_n\}$ such that $g H g^{-1} = H$ for each $g \in S$ then
$$N_G(H) = \bigcup_{g \in S} H g H = \bigcup_{g \in S} g H.$$

\hfill $\square$

\begin{theorem}
\label{comm_ind}
Let $G$ and $H$ be finitely generated non-abelian subgroups of $F^{\mathbb{Z}[t]}$
such that $H \le G$. Then $|\comm{G}{H} : H| < \infty$.
\end{theorem}
{\it Proof.} By Lemma \ref{commensurator} we have $C = \comm{G}{H} = \langle H, g_1,\ldots,g_n\rangle$
for some $g_1,\ldots, g_n \in G$. Consider
$$H' = \langle H \cup H^{g_1} \cup \cdots \cup H^{g_n}\rangle \le G.$$
It is easy to see that $H' \lhd C$, so, $C \le N_G(H')$. Therefore
$|C : H'| < \infty$. Finally, $|H' : H| < \infty$ by Theorem \ref{th:3.1}.

\hfill $\square$

\begin{cor}
\label{cor:comm_ind}
Let $G$ be a finitely generated fully residually free group and let $H$ be its
finitely generated non-abelian subgroup. Then $|\comm{G}{H} : H| < \infty$.
\end{cor}

\begin{cor}
\label{cor:bound_ind}
Let $G$ be a finitely generated fully residually free group and let $H$ be its
finitely generated non-abelian subgroup. Then there exists an effectively computable
natural number $n(H)$ such that for every $K \le G$ containing $H$, if $|K : H| < \infty$
then $|K : H| < n(H)$.
\end{cor}

\end{document}